\documentclass[english]{amsart}

\usepackage{amsmath}
\usepackage{amssymb}
\usepackage{amsthm}
\usepackage{amscd}
\usepackage{babel}
\usepackage[latin1]{inputenc}
\usepackage{graphicx}

\renewcommand{\subjclassname}{AMS \textup{2000} Mathematics Subject
Classification\ }

\newtheorem*{teor}{Theorem}
\newtheorem*{coro}{Corollary}
\newtheorem*{prop}{Proposition}
\newtheorem*{lem}{Lemma}
\theoremstyle{definition}
\newtheorem*{exa}{Example}
\newtheorem*{rem}{Remark}

\author{Antonio M. Oller}
\title{The dying rabbit problem revisited}
\address{Departamento de Matem\'{a}ticas, Universidad de Zaragoza\\
C/Pedro Cerbuna 12, 50009 Zaragoza (Espa\~{n}a)} \email{oller@unizar.es}
\begin{document}
\maketitle

\begin{abstract}
In this paper we study a generalization of the Fibonacci sequence in which rabbits are mortal and take more that two months to become mature. In particular we give a general recurrence relation for these sequences (improving the work in \cite{HOG}) and we calculate explicitly their general term (extending the work in \cite{MIL}). In passing, and as a technical requirement, we also study the behavior of the positive real roots of the characteristic polynomial of the considered sequences.
\end{abstract}
\subjclassname{11B39, 11C08}

\section{Introduction}
Fibonacci numbers arose in the answer to a problem proposed by Leonardo de Pisa who asked for the number of rabbits at the $n^{th}$ month if there is one pair of rabbits at the $0^{th}$ month which becomes mature one month later and that breeds another pair in each of the succeeding months, and if these new pairs breed again in the second month following birth. It can be easily proved by induction that the number of pairs of rabbits at the $n^{th}$ month is given by $f_{n}$, with $f_n$ satisfying the recurrence relation:
$$f_0=f_1=1.$$ $$f_n=f_{n-1}+f_{n-2},\ \forall n\geq2.$$
It is not the point here to state any of the many properties of these numbers (see \cite{VAJ} for a good account of them), nevertheless we will recall that if $r_1<r_2$ are the roots of the polynomial $g(x)=x^2-x-1$ then we can see that: \begin{equation}\label{k2}\displaystyle{f_n=\frac{r_1^n}{r_1-r_2}+\frac{r_2^n}{r_2-r_1}}.\end{equation}

In \cite{MIL} the $k$-generalized Fibonacci numbers $f_n^{(k)}$ are defined as follows:
$$f_n^{(k)}=1,\ \forall 0\leq n\leq k-1.$$ $$f_n^{(k)}=\sum_{i=1}^k f_{n-i}^{(k)},\ \forall n\geq k.$$
In this paper Miles proves, among other results, that if $r_{1},\dots,r_{k}$ are the (distinct) roots of $g_k(x)=x^k-x^{k-1}-\dots-x-1$ then:
\begin{equation}f_n^{(k)}=\sum_{i=1}^k\left(\prod_{\substack{i\neq j\\ 1\leq j\leq k}}(r_i-r_j)^{-1}\right)r_i^n.\end{equation}
which, of course reduces to (\ref{k2}) if we set $k=2$. 

Later, in \cite{HOG}, Hoggat and Lind consider the so called ``dying rabbit problem'', previously introduced in \cite{AL1} and studied in \cite{AL2} or \cite{COH}, which consists in letting rabbits die. In the mentioned paper Hoggat and Lind consider a pair of rabbits at the $0^{th}$ time point which produces $B_n$ pairs at the $n^{th}$ time point and which die at the $m^{th}$ time point after birth. Then if $T_n$ is the total number of live pairs of rabbits at the $n^{th}$ time point, and defining $\displaystyle{B(x)=\sum_{n\geq0}B_nx^n}$, $D(x)=x^m$ and $\displaystyle{T(x)=\sum_{n\geq 0}T_nx^n}$, they prove that
$$T(x)=\frac{1-D(x)}{(1-x)(1-B(x))}.$$
and use this formula to find recurrence relations for $T_n$ in various cases.

The goal of this paper is to give a new look at the dying rabbit problem. In the second section we study a family of polynomials, focusing on the behavior of their positive roots. Although motivated by technical requirements, this study turns out to be of intrinsic interest. In the third section we will find a general recurrence relation for the sequence arising in this problem (which is only given in \cite{HOG} for some particular cases) and we will deduce an explicit formula (which also generalizes the work by Miles) for the total number of live pairs at the $n^{th}$ time point. Finally, in an appendix, we give a procedure written using $\textrm{Maple}^{\circledR}$ to calculate terms of the considered sequences.

\section{A family of polynomials and their roots}
Given natural numbers $h,k\geq1$ we define the following polynomial:
$$g_{k,h}(x)=x^{k+h-1}-x^{k-1}-\dots-x-1.$$
In this section we will study, in some sense, the behavior of the roots of this polynomial in terms of $k$ and $h$. In particular we will be interested in the unique positive real root of $g_{k,h}(x)$.

\begin{prop}
$g_{k,h}(x)$ has a unique positive real root $\alpha_{k,h}$ which lies in $(1,2)$.
\end{prop}
\begin{proof}
Just apply Descartes' rule of signs together with Bolzano's theorem.
\end{proof}

\begin{rem}
Note that $\alpha_{1,h}=1$ for all $h\geq 1$.
\end{rem}

As a consequence of the previous proposition we have the following lemma. The proof is elementary and we omit. 

\begin{lem}
Let $h,k\geq1$ and $0\leq y\in\mathbb{R}$, then $y>\alpha_{k,h}$ if and only if $g_{k,h}(y)>0$.
\end{lem}

Now, for every fixed $h$ we consider the sequence $\{\alpha_{k,h}\}_{k\geq1}$. Similarly, for every fixed $k$ we have a sequence $\{\alpha_{k,h}\}_{h\geq1}$. The following proposition studies the monotony of these sequences.

\begin{prop}
Let $h,k\geq 1$. Then:
\begin{enumerate}
\item $\alpha_{k,h}<\alpha_{k+1,h}$.
\item $\alpha_{k,h}\geq\alpha_{k,h+1}$, and the equality holds if and only if $k=1$.
\end{enumerate}
\end{prop}
\begin{proof}
\begin{enumerate}
\item By definition we know that $g_{k+1,h}(\alpha_{k+1,h})=0$. Now, we have that
$\displaystyle{\alpha_{k+1,h}^{k+h-1}=\frac{\alpha_{k+1,h}^{k+h}}{\alpha_{k+1,h}}=\frac{\alpha_{k+1,h}^k+\alpha_{k+1,h}^{k-1}+\cdots+\alpha_{k+1,h}+1}{\alpha_{k+1,h}}}>\alpha_{k+1,h}^{k-1}+\cdots+\alpha_{k+1,h}+1$,
so $g_{k,h}(\alpha_{k+1,h})>0$ and the result follows from the previous lemma.
\item Again by definition $g_{k,h}(\alpha_{k,h})=0$ and we can write
$g_{k,h+1}(\alpha_{k,h})=\alpha_{k,h}^{k+h}-\alpha_{k,h}^{k-1}-\cdots-\alpha_{k,h}-1=\alpha_{k,h}^{k+h}-\alpha_{k,h}^{k+h-1}=\alpha_{k,h}^{k+h-1}(\alpha_{k,h}-1)\geq0$, with the equality holding if and only if $k=1$. An application of the lemma completes the proof.
\end{enumerate}
\end{proof}

Now, as we have seen that both sequences $\{\alpha_{k,h}\}_{k\geq1}$ and $\{\alpha_{k,h}\}_{h\geq1}$ are monotonic, and as they are clearly bounded, they must be convergent. The next result is devoted to calculate their limits.

\begin{prop}
\begin{enumerate}
\item $\displaystyle{\lim_{k\rightarrow\infty}\alpha_{k,h}=\alpha_h}$ for all $h\geq 1$, where $\alpha_h$ is the unique positive root of the polynomial $p_h(x)=x^h-x^{h-1}-1$.
\item $\displaystyle{\lim_{h\rightarrow\infty}\alpha_{k,h}=1}$ for all $k\geq1$.
\end{enumerate}
\end{prop}
\begin{proof}
\begin{enumerate}
\item Let us fix $h\geq1$. Then for any $k\geq2$ we have $\alpha_{k,h}^{k+h-1}=1+\alpha_{k,h}+\cdots+\alpha_{k,h}^{k-1}=\displaystyle{\frac{\alpha_{k,h}^k-1}{\alpha_{k,h}-1}}$ and thus $\alpha_{k,h}^h-\alpha_{k,h}^{h-1}-1=\displaystyle{\frac{-1}{\alpha_{k,h}^k}}$. Now, as we know that $\displaystyle{\alpha_h=\lim_{k\rightarrow\infty}\alpha_{k,h}>1}$ it is enough to take limits in the previous equality to obtain the result.
\item Let us fix now $k\geq1$. Then for any $h\geq2$ we have $\alpha_{k,h}^{k+h-1}=1+\alpha_{k,h}+\cdots+\alpha_{k,h}^{k-1}$ so, we take logarithms in both sides to obtain the equality  $\displaystyle{\log\alpha_{k,h}=\frac{\log(1+\alpha_{k,h}+\cdots+\alpha_{k,h}^{k-1})}{k+h-1}}$. Finally, if we call $\displaystyle{\beta_k=\lim_{h\rightarrow\infty}\alpha_{k,h}}$ and take limits in the previous expression we arrive at $\log\beta_k=0$ for every $k\geq1$ and the proof is complete.
\end{enumerate}
\end{proof}

The previous propositions can be summarized in the following diagram:
$$\begin{array}{ccccccccccccccc}
\alpha_{1,1}&<&\alpha_{2,1}& < & \alpha_{3,1} & < & \alpha_{4,1}  & < & \dots & < & \alpha_{k,1} & < & \dots & \rightarrow & \alpha_1 \\
\shortparallel & &  \vee &  & \vee &  & \vee   &  &  &  & \vee &  &  &  & \vee\\
   \alpha_{1,2}&<&\alpha_{2,2}& < & \alpha_{3,2} & < & \alpha_{4,2}  & < & \dots & < & \alpha_{k,2} & < & \dots & \rightarrow & \alpha_2 \\
   \shortparallel& & \vee &  & \vee &  & \vee  &  &  &  & \vee &  &  &  & \vee\\
  \alpha_{1,3}& <& \alpha_{2,3}& < & \alpha_{3,3} & < & \alpha_{4,3} & < & \dots & < & \alpha_{k,3} & < & \dots & \rightarrow & \alpha_3 \\
   \shortparallel &  & \vee &  & \vee &  & \vee &  &  &  & \vee& & & & \vee\\
  \alpha_{1,4} &<&\alpha_{2,4}& < & \alpha_{3,4} & < & \alpha_{4,4}  & < & \dots & < & \alpha_{k,4} & < & \dots & \rightarrow & \alpha_4 \\
  \shortparallel &  & \vee &  & \vee &  & \vee &  &  &  & \vee & & & & \vee\\
  \alpha_{1,5} &<&\alpha_{2,5}& < & \alpha_{3,5} & < & \alpha_{4,5}  & < & \dots & < & \alpha_{k,5} & < & \dots & \rightarrow & \alpha_5 \\
   \shortparallel & & \vee&  & \vee &  & \vee &  &  &  & \vee &  &  &  & \vee \\
   \vdots &  & \vdots& &\vdots &  & \vdots &  & \ddots &  & \vdots &  & \ddots &  &\vdots \\
   \shortparallel &  & \vee & &\vee &  & \vee &  &  &  & \vee &  &  &  & \vee \\
  \alpha_{1,h} &<&\alpha_{2,h}& < & \alpha_{3,h} & < & \alpha_{4,h}  & < & \dots & < & \alpha_{k,h} & < & \dots & \rightarrow & \alpha_h \\
  \shortparallel &  & \vee & & \vee &  & \vee  &  &  &  & \vee &  &  &  & \vee\\
   \vdots &  & \vdots & &\vdots &  & \vdots  &  &  &  & \vdots &  &  &  &\vdots \\
   \shortparallel &  & \downarrow & & \downarrow &  & \downarrow  &  &  &  & \downarrow &  &  &  & \downarrow\\
   1 & = & 1 & = & 1 & = & 1 & =&\dots & = & 1 & = & \dots & \rightarrow & 1
\end{array}$$
Where $\alpha_h$ is the unique positive root of $p_h(x)=x^h-x^{h-1}-1$ and every inequality is strict.

Before we go on, we introduce a result by Cauchy (see \cite{CAU}) which will be useful in the forthcoming. This result gives a bound on the modulus of the roots of a polynomial with complex coefficients. We present it without proof.

\begin{teor}
Let $f(z)=z^n+a_{n-1}z^{n-1}+\cdots+a_1z+a_0$ be a polynomial with complex coefficients ($a_i\neq0$ for at least one $i$) and let $Z[f(z)]$ the set of its complex roots. Then $Z[f(z)]\subset\{z\in\mathbb{C}\ |\ |z|\leq\eta\}$ where $\eta$ is the unique positive root of the real polynomial $\widetilde{f}(x)=x^n-|a_{n-1}|x^{n-1}-\cdots-|a_1|x-|a_0|$.
\end{teor}

\begin{coro}
$Z[g_{k,h}(x)]\subset\{z\in\mathbb{C}\ |\ |z|\leq\alpha_{k,h}\}$.
\end{coro}
\begin{proof}
Since $\widetilde{g}_{k,h}(x)=g_{k,h}(x)$, it is enough to apply the previous theorem.
\end{proof}
 
Now, we can refine the previous corollary in the following way.

\begin{prop}
If $g_{k,h}(w)=0$ then $|w|\leq\alpha_{k,h}$ and moreover, $|w|=\alpha_{k,h}$ if and only if $w=\alpha_{k,h}$. In other words, $\alpha_{k,h}$ is the largest root of $g_{k,h}(x)$.
\end{prop}
\begin{proof}
In what follows we will put $\alpha=\alpha_{k,h}$. Let $w$ be a root of $g_{k,h}(x)$ such that $|w|=\alpha$. As $g_{k,h}(w)=0$ we have that $w^{k-1}(w^h-1)=w^{k-2}+\cdots+w+1$ and, in particular, $|w^{k-1}|^2|w^h-1|^2|w-1|^2=|w^{k-1}-1|^2$. If we put $w=\alpha\cos\theta+i\alpha\sin\theta$ we obtain the equality
$\alpha^{2k-2}\left((\alpha^h\cos h\theta-1)^2+\alpha^{2h}\sin^2h\theta\right)\left((\alpha\cos\theta-1)^2+\alpha^2\sin^2\theta\right)=\left((\alpha^{k-1}\cos(k-1)\theta-1)^2+\alpha^{2k-2}\sin^2(k-1)\theta\right)$
and operating we get
$\alpha^{2k-2}(\alpha^{2h}+1-2\alpha^h\cos h\theta)(\alpha^2+1-2\alpha\cos\theta)=\alpha^{2k-2}+1-2\alpha^{k-1}\cos(k-1)\theta$.

Let us now define polynomials $p(x)$ and $q(x)$ by:
$$p(x)=\sum_{i=0}^{\left[\frac{h}{2}\right]}\binom{h}{2i}x^{h-2i}(1-x^2)^{2i},\quad q(x)=\sum_{i=0}^{\left[\frac{k-1}{2}\right]}\binom{k-1}{2i}x^{k-1-2i}(1-x^2)^{2i}$$
then, it is easy to see that $\cos h\theta=p(\cos\theta)$ and $\cos(k-1)\theta=q(\cos\theta)$.
Moreover, $p(1)=q(1)=1$ and if we define $r(x)$ to be
$$r(x)=\alpha^{2k-2}(\alpha^{2h}+1-2\alpha^hp(x))(\alpha^2+1-2\alpha x)-\left(\alpha^{2k-2}+1-2\alpha^{k-1}q(x)\right)$$
then, $g_{k,h}(\alpha)=0$ implies $r(1)=0$. Also we can see that $r(x)$ has no roots in $(-1,1)$. Thus, if $w\in\mathbb{C}$ is a root of $g_{k,h}(x)$ with $|w|=\alpha$ it must be real and as $g_{k,h}(-\alpha)\neq0$, we conclude that $$Z[g_{k,h}(x)]\setminus\{\alpha_{k,h}\}\subset\{z\in\mathbb{C}\ |\ |z|<\alpha_{k,h}\}$$ and the proof is complete.
\end{proof}

We will finish this section with the following proposition which will be of great technical importance in the next section.

\begin{prop}
All the roots of $g_{k,h}$ are distinct.
\end{prop}
\begin{proof}
We will show that $g_{k,h}(x)$ and $g'_{k,h}(x)$ have no common root. To see so it is enough to prove that $\textrm{g.c.d.}(g_{k,h},g'_{k,h})$ is a constant. But we can use Euclid's algorithm repeatedly (multiplying by appropriate constants in every step if necessary) to arrive to the result. 
\end{proof}

\section{The dying rabbit sequence}
As we mentioned in the introduction, we are interested in generalizing the Fibonacci sequence by considering that rabbits become mature $h$ months after their birth and that they die $k$ months after their matureness. We will denote by $C_n^{(k,h)}$ the number of couples of rabbits at the $n^{th}$ month. Obviously we have:
$$C_0^{(k,h)}=\dots=C_{h-1}^{(k,h)}=1.$$

Now let us denote by $C_n^{(h)}$ the recurrence sequence defined by:
$$C_0^{(h)}=\dots=C_{h-1}^{(h)}=1,\ C_n^{(h)}=C_{n-1}^{(h)}+C_{n-h}^{(h)}\ \forall n\geq h$$
then an easy induction argument let us see that:
$$C_n^{(k,h)}=\begin{cases}
C_n^{(h)}, & \text{if $0\leq n \leq k+h-2$} \\
C_{n-h}^{(k,h)}+C_{n-h-1}^{(k,h)}+\cdots+C_{n-k-h+1}^{(k,h)}, & \text{if $n > k+h-2$} \\
\end{cases}$$

\begin{exa}
(See the appendix)
\begin{itemize}
\item If $k=3$ and $h=2$, the beginning terms of $C_n^{(3,2)}$ are:
$$1,1,2,3,4,6,9,13,19,28,41,\dots$$
\item If $k=7$ and $h=4$, then the beginning terms of $C_n^{(7,4)}$ are:
$$1,1,1,1,2,3,4,5,7,10,13,17,23,32,\dots$$
\end{itemize}
\end{exa}

\begin{rem}
The characteristic polynomial of the recurrence sequence $C_n^{(k,h)}$ is precisely the polynomial $g_{k,h}(x)$ studied in the previous section.
\end{rem}

If we denote by $r_1,r_2,\dots,r_{k+h-1}$ the (distinct) complex roots of $g_{k,h}(x)$ it follows from section 2 and from well-known facts from the theory of recurrence sequences that there exist constants $a_1,a_2,\dots,a_{k+h-1}$ such that:
$$C_n^{(k,h)}=a_1r_1^n+a_2r_2^n+\cdots+a_{k+h-1}r_{k+h-1}^n$$
where we can suppose that $r_1=\alpha_{k,h}$. In particular we can calculate those constants solving the system of linear equations given by:
\begin{equation}\label{eq}\sum_{i=1}^{k+h-1}a_ir_i^l=C_l^{(h)},\quad 0\leq l\leq k+h-2.\end{equation}
which can be expressed matricially:
$$\begin{pmatrix}
1 & 1 & \cdots & 1\\ r_1 & r_2 & \cdots & r_{k+h-1}\\ \dots & \vdots & \ddots & \vdots\\ r_1^{k+h-2} & r_2^{k+h-2} & \cdots & r_{k+h-1}^{k+h-2}\end{pmatrix}\begin{pmatrix}a_1\\ a_2\\ \vdots\\ a_{k+h-1}\end{pmatrix}=\begin{pmatrix} C_0^{(h)}\\ C_1^{(h)}\\ \vdots\\ C_{k+h-2}^{(h)}\end{pmatrix}$$
and which has unique solution because all the $r_i$ are distinct.

To solve this system of equations we will use Cramer's rule. Recall that if we put
$$V=\begin{vmatrix}
1 & 1 & \cdots & 1\\ r_1 & r_2 & \cdots & r_{k+h-1}\\ \dots & \vdots & \ddots & \vdots\\ r_1^{k+h-2} & r_2^{k+h-2} & \cdots & r_{k+h-1}^{k+h-2}\end{vmatrix}=\prod_{k+h-1\geq i>j\geq1} (r_i-r_j)$$
$$D_n=\begin{vmatrix}1 & \dots & 1 & C_0 & 1 & \dots & 1\\ r_1 & \dots & r_{n-1} & C_1 & r_{n+1} & \dots & r_{k+h-1}\\ \vdots & \ddots & \vdots & \vdots & \vdots & \ddots & \vdots \\ r_1^{k+h-2} & \dots & r_{n-1}^{k+h-2} & C_{k+h-2} & r_{n+1}^{k+h-1} & \dots & r_{k+h-1}^{k+h-2}\end{vmatrix}$$
then $a_n=\displaystyle{\frac{D_n}{V}}$. So it is enough to find $D_n$ for $n=1,\dots,k+h-1$. We will work out the case $n=1$ completely, the other cases being analogous. Also note that we have replaced the values $C_j^{(h)}$ ($j=0,\dots,k+h-2$) by arbitrary constants $C_j$ ($j=0,\dots,k+h-2$), as their value is of no importance when solving (\ref{eq}) formally.

To calculate $D_1$ we first need the following generalization of Vandermonde determinant which can be found in \cite{ERN} (Lemma 2.1).

\begin{lem}
If $e_n$ is the $n^{th}$ elementary symmetric polynomial in the variables $\{x_1,\dots,x_n\}$, then:
$$\begin{vmatrix} 1 &\dots & 1 & 1\\ x_1 & \dots & x_{n-1} & x_n \\ \vdots & \ddots & \vdots & \vdots\\ \widehat{x_1^l} & \dots & \widehat{x_{n-1}^l} & \widehat{x_n^l}\\ \vdots & \ddots & \vdots & \vdots\\ x_1^n & \dots & x_{n-1}^n & x_n^n\end{vmatrix}=\left(\prod_{n\geq i>j\geq 1}(x_i-x_j)\right)e_{n-l}(x_1,\dots,x_n).$$

\end{lem}

Now, if we apply this lemma and we expand the determinant $D_1$ by its first column we obtain:
$$D_1=\sum_{l=0}^{k+h-2}(-1)^lC_l\left(\prod_{k+h-1\geq i>j\geq 2}(r_i-r_j)\right)e_{k+h-2-l}(r_2,\dots, r_{k+h-1})$$
and, consequently:
\begin{equation}\label{a1}a_1=\frac{D_1}{V}=\frac{1}{\displaystyle{\prod_{k+h-1\geq i\geq2}(r_i-r_1)}}\sum_{l=0}^{k+h-2}(-1)^lC_le_{k+h-2-l}(r_2,\dots, r_{k+h-1}).\end{equation}

So, we are now interested in calculating the values $e_j(r_2,\dots,r_{k+h-1})$ for $0\leq j\leq k+h-2$. From Cardano's formulae and considering that $r_1,\dots,r_{k+h-1}$ are the roots of $g_{k,h}(x)$ we have that:
$$e_0(r_1,\dots,r_{k+h-1})=1.$$
$$e_1(r_1,\dots,r_{k+h-1})=\dots=e_{h-1}(r_1,\dots,r_{k+h-1})=0.$$
$$e_s(r_1,\dots,r_{k+h-1})=(-1)^{s+1}\quad \forall h\leq s\leq k+h-1.$$

On the other hand, the following lemma is easy to prove.

\begin{lem}
$\displaystyle{e_t(x_2,\dots,x_n)=\sum_{i=1}^{n-t}(-1)^{i+1}\frac{e_{t+i}(x_1,\dots,x_n)}{x_1^i}\quad \forall 0\leq t<n}$.
\end{lem}

We put this together to obtain:
$$e_0(r_2,\dots,r_{k+h-1})=1.$$
$$e_s(r_2,\dots,r_{k+h-1})=(-1)^{s}\sum_{i=1}^k\frac{1}{r_1^{i+h-1-s}}\quad \forall 1\leq s\leq h-1.$$
$$e_s(r_2,\dots,r_{k+h-1})=(-1)^{s}\sum_{i=1}^{k+h-1-s}\frac{1}{r_1^i}\quad \forall h\leq s\leq k+h-2.$$
and summing the geometric series:
$$e_0(r_2,\dots,r_{k+h-1})=1.$$
$$e_s(r_2,\dots,r_{k+h-1})=(-1)^{s}\frac{r_1^k-1}{r_1^{k+h-1-s}(r_1-1)}\quad \forall 1\leq s\leq h-1.$$
$$e_s(r_2,\dots,r_{k+h-1})=(-1)^{s}\frac{r_1^{k+h-1-s}-1}{r_1^{k+h-1-s}(r_1-1)}\quad \forall h\leq s\leq k+h-2.$$
Finally if we substitute in (\ref{a1}) we get: 
$$a_1=\frac{(-1)^{k+h}}{\displaystyle{\prod_{k+h-1\geq i>2}(r_i-r_1)}}\left[\sum_{l=0}^{k-2}C_l\frac{r_1^{l+1}-1}{r_1^{l+1}(r_1-1)}+\sum_{l=k-1}^{k+h-3}C_l\frac{r_1^k-1}{r_1^{l+1}(r_1-1)}+C_{k+h-2}\right].$$

Reasoning in a similar way and taking into account the symmetry of the $e_s$ we can calculate $a_n$ for every $1\leq n\leq k+h-1$. In fact:
$$a_n=\frac{(-1)^{k+h+n-1}}{\displaystyle{\prod_{i>n}(r_i-r_n)}\prod_{n>j}(r_n-r_j)}\left[\sum_{l=0}^{k-2}C_l\frac{r_n^{l+1}-1}{r_n^{l+1}(r_n-1)}+\sum_{l=k-1}^{k+h-3}C_l\frac{r_n^k-1}{r_n^{l+1}(r_n-1)}+C_{k+h-2}\right].$$

\begin{rem}
Observe that the previous considerations are only valid for $k,h\geq2$. For the case $h=1$ we refer to Miles' paper \cite{MIL} and the case $k=1$ is trivial.
\end{rem}

\begin{rem}
It is interesting to observe that $a_1\neq 0$. As a consequence and recalling that $|r_i|<|r_1|$ for all $i\geq2$ we have that $\displaystyle{\lim_{n\rightarrow\infty}\frac{C_{n+1}^{(k,h)}}{C_n^{(k,h)}}=r_1=\alpha_{k,h}}$. This generalizes
 the fact that $\displaystyle{\frac{f_{n+1}}{f_n}=\Phi}$ where $f_n$ is the $n^{th}$ Fibonacci number and $\Phi$ is de golden section (note that $\alpha_{2,1}=\Phi$). The previous expression for $a_n$ also generalizes the work by Miles just by setting $h=1$.
\end{rem}

\begin{exa}
(Padovan sequence)\newline
Recall that the so-called Padovan sequence is defined by
$$P_0=P_1=P_2=1.$$ $$P_n=P_{n-2}+P_{n-3},\quad \forall n\geq3.$$
Thus, it is clear that in our notation $P_n=C_n^{(2,2)}$ with the initial conditions $C_0=C_1=C_2=1$. So we can apply our previous results to obtain that:
$$P_n=\frac{r_1^2+r_1+1}{2r_1+3}r_1^n+\frac{r_2^2+r_2+1}{2r_2+3}r_2^n+\frac{r_3^2+r_3+1}{2r_3+3}r_3^n.$$
which was already known to hold.\newline
If we keep the same recurrence relation but replace the initial conditions by $P_0=3,P_1=0,P_2=2$ we obtain the so-called Perrin sequence, whose general term can be again calculated with our formulas to obtain:
$$P_n=r_1^n+r_2^n+r_3^n.$$
Finally, if we keep our original initial conditions, that is, $C_0^{(2,2)}=C_1^{(2,2)}=1, C_2^{(2,2)}=2$, then the general term of our Padovan-Perrin like sequence turns out to be:
$$C_n^{(2,2)}=\frac{(r_1+1)^2}{2r_1+3}r_1^n+\frac{(r_2+1)^2}{2r_2+3}r_2^{n}+\frac{(r_3+1)^2}{2r_3+3}r_3^{n}.$$
\end{exa}

\section{Appendix}
In this appendix we give a short and easy procedure, written with $\textrm{Maple}^{\circledR}$, which calculates any number of terms of $C_n^{(k,h)}$. It goes as follows:

\begin{verbatim}
dr:=proc(k,h,t)
local i;
for i from 0 by 1 to h-1 do c(i):=1 end do;
for i from h by 1 to k+h-2 do c(i):=c(i-1)+c(i-h); end do;
for i from k+h-1 to t do c(i):=sum(c(n), n=i-k-h+1..i-h); end do;
print(seq(c(n),n=0..t));
end proc:
\end{verbatim}

\end{document}